\documentclass[10pt]{article}

\usepackage{amssymb,amsmath,graphics}

\numberwithin{equation}{section}

\allowdisplaybreaks

\newenvironment{Proof}{\removelastskip\par\medskip
\noindent{\em Proof.}
\rm}{\penalty-20\null\hfill$\square$\par\medbreak}

\def\real{{\mathord{\mathbb R}}}

\newtheorem{prop}{Proposition}

\def\Dom{{\mathrm{{\rm Dom ~}}}}

\def\deltaa{{\mathord{\mathrm {Delta}}}}
\def\rhoa{{\mathord{{\mathrm {Rho}}}}}
\def\vegaa{{\mathord{{\mathrm {Vega}}}}}
\def\gammaa{{\mathord{{\mathrm {Gamma}}}}}

\begin{document}
 ~

\begin{center}
{
\Large
\baselineskip=0.9cm
\textbf{Computations of Greeks in stochastic volatility models via
the Malliavin calculus}} \normalsize

\end{center}
\
\baselineskip=0.5cm
\begin{center}
 \large Youssef El-Khatib, {\em
Department of Mathematical Sciences,Youssef\_Elkhatib@uaeu.ac.ae\\
\vspace{.2cm}
U.A.E University, Al-Ain, P.O.Box 17555, U.A.E}
\normalsize
\end{center}
\begin{abstract}
\baselineskip0.5cm We compute Greeks for stochastic volatility models
driven by Brownian informations. We use the Malliavin method introduced by \cite{fournie} for
deterministic volatility models.
\end{abstract}
\baselineskip=0.5cm \noindent{\bf Keywords: Stochastic volatility, Greeks, Malliavin
calculus}

\baselineskip0.7cm
\section{INTRODUCTION}
\label{intro} The application of the Malliavin calculus to the
computations of price sensitivities were introduced by
\cite{fournie} for models with deterministic volatility. In this
work, we compute the Greeks, for stochastic volatility models where
the underlying asset price is driven by Brownian information. We
consider stochastic volatility models, since these models, unlike
those with deterministic volatility, take into account the
{\it{smile}} effect.\\ Let $(B_t)_{t\in [0,T]}$ and $(B^{'}_t)_{t\in
[0,T]}$ be two independent Brownian motions. We work in a filtered
probability space $(\Omega, {\mathcal{F}},({\mathcal{F}}_t)_{t\in
[0,T]}, P)$, where $({\mathcal{F}}_t)_{t\in [0,T]}$ is the naturel
filtration generated by $B$ and $B^{'}$. We consider a market with
two assets, a riskless one with price $(e^{\int_0^t r_s ds})_{t\in
[0,T]}$, where $r$ is deterministic and denotes the interest rate.
And a risky asset $(S_t)_{t\in [0,T]}$ to which is related an option
with payoff $f(S_T)$. $(S_t)_{t\in [0,T]}$ has a stochastic
volatility and is given by
$$
\frac{dS_t}{S_t}=\mu_t dt+\sigma(t,Y_t)dB_t,\ \ \ \mbox{where}\ \ \
dY_t=\mu^Y_t dt+\sigma^Y_t [\rho dB_t+dB^{'}_t ],\ \ \ t \in[0,T],
$$
with $S_0=x>0$ and $Y_0=y \in \real$. $\mu$, $\mu^{Y}$, $\sigma^Y$
are deterministic functions, $\rho \in \real$ and $\sigma \in
{\mathcal C}^2([0,T]\times \real)$ such that for any $t\in [0,T]$,
$\sigma(t,.)\neq 0$. The market considered here is incomplete. There
is an infinity of E.M.M -Equivalent Martingale Measure- (i.e a
probability equivalent to $P$ under which the actualized price $(S_t
e^{-\int_0^t r_s ds})_{t \in [0,T]}$ is a martingale). Let $Q$ be a fixed
$P$-E.M.M. $Q$ is identified by its Radon-Nikodym density w.r.t $P$, denoted
$\rho_T$ and given by
$$ \rho_T=\exp\left(\int_0^T \alpha_s dB_s+\beta_s dB^{'}_s
-\frac{1}{2}\int_0^T (\alpha^2_s+\beta^2_s)ds\right),$$ where
$(\alpha_{t})_{t\in [0,T]}$ and $(\beta_{t})_{t\in [0,T]}$ are two
predictable processes s.t. $\alpha_t= -\frac{\mu_t-r_t}{\sigma(t,Y_t)}$, and $\beta$
is arbitrary. Now let, for any $t\in
 [0,T]$, $W_t=B_t -\int_0^t \alpha_s ds$ and
$W^{'}_t=B^{'}_t -\int_0^t \beta_s ds$ then by the Girsanov theorem
$W$ and $W^{'}$ are two $Q$-Brownian motions. In the following we
work with a fixed $P$-E.M.M. $Q$ and we will use $E[.]$ (instead of
$E_{Q}$[.]) as the expectation under the probability $Q$. We have,
under $Q$, for any $t\in [0,T]$
\begin{equation}
\label{St}
 \frac{dS_t}{S_t}=r_t dt+\sigma(t,Y_t)dW_t,\ \ \ dY_t=\left(\mu^Y_t +\sigma^Y_t\frac{r_t-\mu_t}{\sigma(t,Y_t)}+\beta_t \sigma^Y_t\right)dt+ \rho \sigma^Y_t
dW_t+\sigma^Y_t dW^{'}_t.
\end{equation}
\section{MALLIAVIN DERIVATIVE ON WIENER SPACE}
In this section we give an introduction to the malliavin derivative
in Wiener space and to its adjoint the Skorohod integral. We refer
for example to \cite{oksendal}. Let $(D^W_t)_{t\in [0,T]}$ be the
Malliavin derivative on the direction of $W$. We denote by
$\mathbb{P}$ the set of random variables $F:\Omega \rightarrow
\real$, such that F has the representation
$$F(\omega)=f\left(\int_0^T f_1(t)dW_t,
  \ldots,\int_0^T f_n(t)dW_t\right),$$ where
  $f(x_1,\ldots,x_n)=\sum_{\alpha}a_{\alpha}x^{\alpha}$ is
  a  polynomial in $n$ variables $x_1,\ldots,x_n$ and
  deterministic functions  $f_i \in L^2([0,T])$ .
Let $\|.\|_{1,2}$ be the norm
$$\|F\|_{1,2}:=\|F\|_{L^2(\Omega)}+\|D_{\cdot}^{W} F\|_{L^2([0,T]\times\Omega)}, \ \ \ F \in L^2(\Omega).$$
Thus the domaine of the operator $D^W$, $\Dom(D^W)$, coincide with
$\mathbb{P}$ w.r.t the norm $\|.\|_{1,2}$. The next proposition will
be useful.
\begin{prop}
\label{p1} Given $F=f\left(\int_0^T f_1(t)dW_t, \ldots,\int_0^T
f_n(t)dW_t\right)\in \mathbb{P}$. We have
$$
D^W_t F =\sum_{k=0}^{k=n}\frac{\partial f}{\partial
x_k}\left(\int_0^T f_1(t)dW_t,\ldots,\int_0^T
f_n(t)dW_t\right)f_k(t).
$$
\end{prop}
To calculate the Mallaivin derivative for It\^o integral, we will
use the following Proposition.
\begin{prop} \label{derivint} Let $(u_t)_{t\in [0,T]}$ be a
${\mathcal{F}}_t-$adapted process, such that $u_t \in \Dom(D^{W})$,
we have
$$D^{W}_t \int_0^T u_r dW_r=\int_t^T (D^{W}_t u_s)dW_s+ u_t.$$
\end{prop}
From now on, for any stochastic process $u$ and for $F\in \Dom(D^W)$
such that $u_. D^W_. F \in L^2([0,T])$  we let
$$D^W_u F:=\langle D^W F,u\rangle_{L^2([0,T])}:=\int_0^T u_t D^W_t
F dt.$$ Let $\delta^W$ be the Skorohod integral in Wiener space. We
have $\delta^W$ is the adjoint of $D^W$ as showing in the next
proposition, moreover its an extension of the It\^o integral
\begin{prop}
a) Let $u \in \Dom(\delta^W)$ and $F\in \Dom(D^W)$, we have $E[D^W_u
F]\leq C(u)\|F\|_{1,2}$, and $E[F\delta^W(u)]=E[D^W_u F]$.\\
b) Consider a $L^2(\Omega \times [0,T])$-adapted stochastic process
$u=(u_t)_{t\in [0,T]}$. We have $\delta^W(u)=\int_0^T u_t dW_t.$\\
c) Let $F\in \Dom(D^W)$ and $u\in \Dom(\delta^W)$ such that $uF \in
\Dom(\delta^W)$ thus $\delta^W(uF)=F\delta^W(u)-D^W_u F.$
\end{prop}
\section{COMPUTATIONS OF THE GREEKS}
The computations of Greeks by Malliavin approach rest on a known
integration by parts formula -cf. \cite{fournie}- given in the
following proposition.
\begin{prop}
\label{mc} Let $I$ be an open interval of $\real$. Let
$(F^\zeta)_{\zeta \in I}$ and $(H^\zeta)_{\zeta \in I}$,
 be two families of random functionals, continuously differentiable
 in $\Dom (D^W)$ in the parameter $\zeta \in I$.
 Let $(u_t)_{t\in [0,T]}$ be a process satisfying
$$D_u^W F^\zeta \neq 0, \quad a.s. \mbox{ on }
 \{ \partial_\zeta F^\zeta \neq 0 \}, \quad \zeta \in I,
$$
 and such that $ u H^\zeta \partial_\zeta F^\zeta /D_u^W F^\zeta$
 is continuous in $\zeta$ in $\Dom (\delta^W)$.
 We have
$$
\frac{\partial}{\partial \zeta}E\left[
 H^\zeta f\left(F^\zeta \right)
 \right]
= E\left[f\left(F^\zeta \right)\left(\frac{H^\zeta \partial_\zeta
F^\zeta}
 {D_u^W F^\zeta }\delta^W(u)-D^W_u
 \left(\frac{H^\zeta \partial_\zeta F^\zeta}{D_u^WF^\zeta
 }\right)+\partial_\zeta H^\zeta\right)\right]
$$
 for any function $f$ such that
 $f\left(F^\zeta \right) \in L^2(\Omega )$,
 $\zeta\in I$.
\end{prop}
Our aim is to compute the Greeks for options with payoff $f(S_T)$,
where $(S_t)_{t\in [0,T]}$ denotes the underlying asset price given
by (\ref{St}). We have
\begin{equation}
\label{Stexp}
 S_{T} = x\exp \left( \int_0^T \sigma(s,Y_s)dW_s
 + \int_0^T (r_s -
 \frac{1}{2}  \sigma^2(s,Y_s) ds) \right),
\end{equation}
Let $\zeta$ be a parameter taking the values: the initial asset
price $x=S_0$, the volatility $\sigma$, or the interest rate $r$.
Let $C=E[f(S^{\zeta}_t)]$ be the price of the option. We will
compute the following Greeks:
\begin{eqnarray*}
& & \deltaa = \frac{\partial C}{\partial x},
  \quad \quad
 \gammaa =\frac{\partial^2 C}{\partial x^2},
  \quad \quad \rhoa=\frac{\partial C}{\partial r},
 \quad \quad
\vegaa = \frac{\partial C}{\partial \sigma},
\end{eqnarray*}
The next Proposition gives the first, second and third order
derivatives of $S_T$ w.r.t $D^{W}$, needed for the computations of
the different Greeks.
\begin{prop}
\label{der}
For $0\leq t\leq T$, we let $G(t,T):=\sigma(t,Y_t)+\int_t^T \frac{\partial \sigma}{\partial
y}(v,Y_v)D^{W}_t Y_v (d{W}_v-\sigma(v,Y_v)dv)$. We have $D^{W}_t S_T=S_T G(t,T)$
and thus
\begin{eqnarray}
\label{der1} D^{W}_u S_T&=&S_T \int_0^T u_t G(t,T)dt\\
\label{der2} D^{W}_u D^{W}_u S_T&=&S_T\left(\left(\int_0^T u_t G(t,T)dt\right)^2+\int_0^T\int_s^Tu_s u_tD_sG(t,T)dtds\right)\\
\nonumber D^W_u D^W_u D^{W}_u S_T &=& S_T\left(\left(\int_0^T u_t G(t,T)dt\right)^3+3\int_0^T u_t G(t,T)dt \int_0^T\int_s^Tu_su_tD_sG(t,T)dtds\right.\\
&&
\label{der3}
\left.+\int_0^T\int_r^T\int_s^T u_ru_su_t D_r D_sG(t,T)dtdsdr\right)
\end{eqnarray}
where
\begin{eqnarray}
\nonumber
 D^W_s G(t,T)&=&\frac{\partial \sigma}{\partial
y}(t,Y_t)D^{W}_s Y_t+\int_t^T D^W_s Y_v D^W_t Y_v
 \frac{\partial}{\partial y}\left[\frac{\partial \sigma}{\partial
y}(v,Y_v)(d{W}_v-\sigma(v,Y_v)dv)\right]\\
\label{derG1} &&+\int_t^T \frac{\partial \sigma}{\partial
y}(v,Y_v)D^{W}_s D^{W}_t Y_v(d{W}_v-\sigma(v,Y_v)dv)\\
\nonumber
D^W_r D^W_s G(t,T)&=&\int_t^T \left((D^W_rD^W_s Y_v D^W_t Y_v+D^W_s Y_v D^W_rD^W_t Y_v)
 \frac{\partial^2 \sigma}{\partial y^2}(v,Y_v)\right.\\
 &&\left.+D^W_r Y_v D^W_s Y_v D^W_t Y_v\frac{\partial^3\sigma}{\partial y^3}(v,Y_v)\right.\\
 \nonumber
 &&\left.+\frac{\partial \sigma}{\partial
y}(v,Y_v)D^W_rD^{W}_s D^{W}_t Y_v+ \frac{\partial^2 \sigma}{\partial
y^2}(v,Y_v)D^W_rY_vD^{W}_s D^{W}_t Y_v\right)d{W}_v\\
\nonumber
&&-\int_t^T \left[(D^W_rD^W_s Y_v D^W_t Y_v+D^W_s Y_v D^W_rD^W_t Y_v)\frac{\partial}{\partial
y}\left(\sigma(v,Y_v)\frac{\partial \sigma}{\partial
y}(v,Y_v)\right)\right.\\
\nonumber
&&+\left.\frac{\partial^2}{\partial
y^2}\left(\sigma(v,Y_v)\frac{\partial \sigma}{\partial
y}(v,Y_v)\right)D^W_rY_vD^{W}_s D^{W}_t Y_v\right.\\
\nonumber
&&\left.+\frac{\partial}{\partial
y}\left(\sigma(v,Y_v)\frac{\partial \sigma}{\partial
y}(v,Y_v)\right)D^{W}_r Y_vD^{W}_s D^{W}_t Y_v\right.\\
&&
\nonumber
+\left.\sigma(v,Y_v)\frac{\partial \sigma}{\partial
y}(v,Y_v)D^{W}_rD^{W}_s D^{W}_t Y_v\right]dv\\
&&
\label{derG2}
+\frac{\partial^2 \sigma}{\partial y^2}(t,Y_t)D^{W}_r Y_t D^{W}_s Y_t+\frac{\partial
\sigma}{\partial y}(t,Y_t)D^{W}_rD^{W}_s Y_t.
\end{eqnarray}
\end{prop}
\begin{Proof}
By the chain rule of $D^{W}_t$ and thanks to
Proposition~\ref{derivint} we obtain
\begin{eqnarray*}
D^{W}_t S_T&=&S_T \left(D^{W}_t \int_0^T \left(r_s - \frac{\sigma^2
(s,Y_s)}{2}\right) ds +D^{W}_t
 \int_0^T \sigma(s,Y_s)dW_s\right)\\
 &=&S_T \left(-\int_t^T D^{W}_t
  \frac{\sigma^2 (s,Y_s)}{2}ds +
 \int_t^T D^{W}_t \sigma(s,Y_s)dW_s
 +\sigma(t,Y_t)\right)=S_T G(t,T),
\end{eqnarray*}
which gives (\ref{der1}). (\ref{der2}) and (\ref{der3}) are
immediate by the chain rule of $D^{W}$. Concerning
(\ref{derG1}) and (\ref{derG2}) we have for $0\leq r \leq s\leq t\leq T$
\begin{eqnarray*}
D^W_s G(t,T)&=&\frac{\partial \sigma}{\partial
y}(t,Y_t)D^{W}_s Y_t+ \int_t^T D^W_s \left(\frac{\partial
\sigma}{\partial y}(v,Y_v)D^{W}_t Y_v \right)d{W}_v\\
&&-\int_t^T D^W_s\left(
 \sigma(v,Y_v)\frac{\partial \sigma}{\partial y}(v,Y_v)D^{W}_t
Y_v\right)dv\\
 &=&\frac{\partial \sigma}{\partial
y}(t,Y_t)D^{W}_s Y_t+\int_t^T \left(D^W_s Y_v D^W_t Y_v
 \frac{\partial^2
\sigma}{\partial y^2}(v,Y_v)+\frac{\partial \sigma}{\partial
y}(v,Y_v)D^{W}_s D^{W}_t Y_v \right)d{W}_v\\
&&-\int_t^T \left[D^W_s Y_v D^W_t Y_v\frac{\partial}{\partial
y}\left(\sigma(v,Y_v)\frac{\partial \sigma}{\partial
y}(v,Y_v)\right)+ \sigma(v,Y_v)\frac{\partial \sigma}{\partial
y}(v,Y_v)D^{W}_s D^{W}_t Y_v\right]dv.
\end{eqnarray*}
And
\begin{eqnarray*}
D^W_r D^W_s G(t,T)&=&\int_t^T D^W_r \left(D^W_s Y_v D^W_t Y_v
 \frac{\partial^2 \sigma}{\partial y^2}(v,Y_v)+\frac{\partial \sigma}{\partial
y}(v,Y_v)D^{W}_s D^{W}_t Y_v \right)d{W}_v\\
&&-\int_t^T D^W_r\left[D^W_s Y_v D^W_t Y_v\frac{\partial}{\partial
y}\left(\sigma(v,Y_v)\frac{\partial \sigma}{\partial
y}(v,Y_v)\right)\right.\\
&&\left.+\sigma(v,Y_v)\frac{\partial \sigma}{\partial
y}(v,Y_v)D^{W}_s D^{W}_t Y_v\right]dv\\
&&+\frac{\partial^2 \sigma}{\partial y^2}(t,Y_t)D^{W}_r Y_t D^{W}_s Y_t+\frac{\partial
\sigma}{\partial y}(t,Y_t)D^{W}_rD^{W}_s Y_t\\
&=&\int_t^T \left((D^W_rD^W_s Y_v D^W_t Y_v+D^W_s Y_v D^W_rD^W_t Y_v)
 \frac{\partial^2 \sigma}{\partial y^2}(v,Y_v)\right.\\
 &&\left.+D^W_r Y_v D^W_s Y_v D^W_t Y_v\frac{\partial^3\sigma}{\partial y^3}(v,Y_v)\right.\\
 &&\left.+\frac{\partial \sigma}{\partial
y}(v,Y_v)D^W_rD^{W}_s D^{W}_t Y_v+ \frac{\partial^2 \sigma}{\partial
y^2}(v,Y_v)D^W_rY_vD^{W}_s D^{W}_t Y_v\right)d{W}_v\\
&&-\int_t^T \left[(D^W_rD^W_s Y_v D^W_t Y_v+D^W_s Y_v D^W_rD^W_t Y_v)\frac{\partial}{\partial
y}\left(\sigma(v,Y_v)\frac{\partial \sigma}{\partial
y}(v,Y_v)\right)\right.\\
&&+\left.\frac{\partial^2}{\partial
y^2}\left(\sigma(v,Y_v)\frac{\partial \sigma}{\partial
y}(v,Y_v)\right)D^W_rY_vD^{W}_s D^{W}_t Y_v\right.\\
&&\left.+\frac{\partial}{\partial
y}\left(\sigma(v,Y_v)\frac{\partial \sigma}{\partial
y}(v,Y_v)\right)D^{W}_r Y_vD^{W}_s D^{W}_t Y_v\right.\\
&&+\left.\sigma(v,Y_v)\frac{\partial \sigma}{\partial
y}(v,Y_v)D^{W}_rD^{W}_s D^{W}_t Y_v\right]dv\\
&&+\frac{\partial^2 \sigma}{\partial y^2}(t,Y_t)D^{W}_r Y_t D^{W}_s Y_t+\frac{\partial
\sigma}{\partial y}(t,Y_t)D^{W}_rD^{W}_s Y_t.
\end{eqnarray*}
\end{Proof}
Concerning $D^W_t Y_v$, it can be explicitly computed when $\beta_v=\beta(v,Y_v)$. We have for $0\leq
t\leq v\leq T$
\begin{eqnarray*}
D^{W}_t Y_v &=&D^{W}_t\left(\int_0^v \left(\mu^Y_\alpha +\sigma^Y_\alpha
\frac{r_\alpha-\mu_\alpha}{\sigma(\alpha,Y_\alpha)}+\beta(\alpha,Y_\alpha) \sigma^Y_\alpha\right)d\alpha+\int_0^v \rho \sigma^Y_\alpha
dW_\alpha+\int_0^v \sigma^Y_\alpha dW^{'}_\alpha\right)
\\
&=&\int_t^v D^{W}_t \left(\mu^Y_\alpha +\sigma^Y_\alpha\frac{r_\alpha-\mu_\alpha}{\sigma(\alpha,Y_\alpha)}+\beta(\alpha,Y_\alpha) \sigma^Y_\alpha\right)d\alpha+
 D^{W}_t \int_0^v \rho \sigma^Y_\alpha
dW_\alpha\\
&=&\int_t^v \sigma^Y_\alpha D^{W}_t \left(
\frac{r_\alpha-\mu_\alpha}{\sigma(\alpha,Y_\alpha)}+\beta(\alpha,Y_\alpha)\right)d\alpha+
  \int_t^v D^{W}_t \rho \sigma^Y_\alpha
dW_\alpha +\rho \sigma^Y_t\\
&=&\rho\sigma^{Y}_t -\int_t^v \sigma^{Y}_\alpha\left(
\frac{r_\alpha-\mu_\alpha}{\sigma^2(\alpha,Y_\alpha)}-\frac{\partial \beta(\alpha,Y_\alpha)}{\partial y}\right)D^{W}_t Y_\alpha d\alpha.
\end{eqnarray*}
So for $t$ fixed in $[0,T]$, the Malliavin derivative of $Y_v$ for $v\in [t,T]$ :
$(D^{W}_t Y_v)_{v\in [t,T]}$, satisfies a stochastic differential
equation, its solution is precisely
$$ D^{W}_t Y_v=\rho \sigma^{Y}_t
\exp\left(-\int_t^v \left(\sigma^{Y}_\alpha
\frac{r_\alpha-\mu_\alpha}{\sigma^2(\alpha,Y_\alpha)}-\frac{\partial \beta(\alpha,Y_\alpha)}{\partial y}\right)d\alpha\right)\ \ \ v\in[t,T].
$$
In this case for $0\leq r\leq s \leq t \leq v \leq T$, $D^W_s D^W_t Y_v$ and $D^W_r D^W_s D^W_t Y_v$ are as follow
\begin{eqnarray*}
D^W_s D^W_t Y_v&=&D^W_t Y_v D^W_s\left(-\int_t^v \left(\sigma^{Y}_\alpha
\frac{r_\alpha-\mu_\alpha}{\sigma^2(\alpha,Y_\alpha)}-\frac{\partial \beta(\alpha,Y_\alpha)}{\partial y}\right)d\alpha\right)\\
&=&D^W_t Y_v  \int_t^v \left(2\sigma^{Y}_\alpha
\frac{r_\alpha-\mu_\alpha}{\sigma^3(\alpha,Y_\alpha)}+\frac{\partial^2 \beta(\alpha,Y_\alpha)}{\partial y^2}\right)D^W_s Y_\alpha d\alpha,\\
D^W_r D^W_s D^W_t Y_s&=&D^W_r \left(D^W_t Y_v  \int_t^v \left(2\sigma^{Y}_\alpha
\frac{r_\alpha-\mu_\alpha}{\sigma^3(\alpha,Y_\alpha)}+\frac{\partial^2 \beta(\alpha,Y_\alpha)}{\partial y^2}\right)D^W_s Y_\alpha d\alpha \right)\\
&=&(D^W_rD^W_t Y_v) \int_t^v \left(2\sigma^{Y}_\alpha
\frac{r_\alpha-\mu_\alpha}{\sigma^3(_\alpha,Y_\alpha)}+\frac{\partial^2 \beta(_\alpha,Y_\alpha)}{\partial y^2}\right)D^W_s Y_\alpha dv_\alpha\\
&&-D^W_t Y_v \int_t^v \left(6\sigma^{Y}_\alpha
\frac{r_\alpha-\mu_\alpha}{\sigma^4(\alpha,Y_\alpha)}-\frac{\partial^3 \beta(\alpha,Y_\alpha)}{\partial y^3}\right)D^W_r Y_\alpha D^W_s Y_\alpha d\alpha\\
&&+D^W_t Y_v \int_t^v \left(2\sigma^{Y}_\alpha
\frac{r_\alpha-\mu_\alpha}{\sigma^3(\alpha,Y_\alpha)}+\frac{\partial^2 \beta(\alpha,Y_\alpha)}{\partial y^2}\right)D^W_rD^W_s Y_\alpha d\alpha.
\end{eqnarray*}
\subsection{$\deltaa$, $\rhoa$, $\vegaa$}
The $\deltaa$, $\rhoa$, $\vegaa$ can be computed using a first order
derivative of $C=E[f(S_T^{\zeta})]$ w.r.t $\zeta$. We have using
Proposition~\ref{mc}
$$
\frac{\partial}{\partial \zeta}E\left[
 f(S_T^\zeta)
 \right]
= E\left[f(S_T^\zeta)\left(\frac{\partial_\zeta S_T^\zeta}
 {D_u^W S_T^\zeta }\delta^W(u)-D^W_u
 \left(\frac{S_T^\zeta \partial_\zeta S_T^\zeta}{D_u^W S_T^\zeta
 }\right)\right)\right].
$$
Next we compute the $\deltaa$, the $\rhoa$ and $\vegaa$ can be computed by the
same way. The $\deltaa$ corresponds to $\zeta=x$, so
$\partial_\zeta S_T=\partial_x S_T=\frac{1}{x}S_T$ and we have
\begin{eqnarray*}
\nonumber
\deltaa&=&E\left[f(S_T)\left(\frac{\partial_x S_T}
 {D_u^W S_T}\delta^W(u)-D^W_u
 \left(\frac{S_T\partial_x S_T}{D_u^W S_T
 }\right)\right)\right]\\
 \nonumber
 &=&E\left[f(S_T)\left(\frac{S_T}
 {x D_u^W S_T}\delta^W(u)-D^W_u
 \left(\frac{S_T^2}{xD_u^W S_T
 }\right)\right)\right]\\
 \nonumber
 &=&\frac{1}{x}E\left[f(S_T)\left(\frac{1}
 {\int_0^T u_t G(t,T)dt}\delta^W(u)-2S_T+
 \frac{S_T^2D_u^WD_u^WS_T}{(D_u^W S_T)^2
 }\right)\right]\\
 \nonumber
 &=&\frac{1}{x}E\left[f(S_T)\left(\frac{1}
 {\int_0^T u_t G(t,T)dt}\delta^W(u)-2S_T\right.\right.\\
 &&\left.\left.+
 \frac{S_T(\int_0^T u_t G(t,T)dt)^2+S_T\int_0^T u_s D_s(\int_0^T u_t G(t,T)dt)ds}{(\int_0^T u_t G(t,T)dt)^2
 }\right)\right]\\
 \label{del}
 &=&\frac{1}{x}E\left[f(S_T)\left(\frac{1}
 {\int_0^T u_t G(t,T)dt}\delta^W(u)-S_T\left(1-
 \frac{\int_0^T u_s (\int_s^T u_t D_sG(t,T)dt)ds}{(\int_0^T u_t G(t,T)dt)^2
 }\right)\right)\right],
\end{eqnarray*}
where $D^W_s G(t,T)$ is given by $(\ref{derG1})$.
\subsection{$\gammaa$}
The $\gammaa$ is computed using the second order derivative of
$C=E[f(S_T)]$ w.r.t $x$ given by
\begin{eqnarray*}
 \label{gam}
\lefteqn{\frac{\partial^2}{\partial x^2}
 E\left[f(S_T) \right]=\frac{\partial}{\partial x}\deltaa=\frac{1}{x}\frac{\partial}{\partial x}E\left[f(S_T)H\right]}\\
 \nonumber
 &=&\frac{1}{x}E\left[f(S_T)\left(\frac{H\partial_x S_T}
 {D_u^W S_T}\delta^W(u)-D^W_u\left(\frac{H\partial_x S_T}{D_u^W S_T}\right)+\partial_x H \right)\right]\\
\nonumber
 &=&\frac{1}{x}E\left[f(S_T)\left(\frac{HS_T}
 {xD_u^W S_T}\delta^W(u)-\left(\frac{H(D^W_uS_T)^2+S_TD^W_uS_TD^W_uH+HS_TD_u^W D_u^W S_T}{x(D_u^W S_T)^2}\right)\right.\right.\\
 &&\left.\left.+\partial_x H \right)\right],
 \end{eqnarray*}
 where
\begin{eqnarray*}
H&=&\frac{1}{\int_0^T u_t G(t,T)dt}\delta^W(u)-S_T\left(1-
 \frac{\int_0^T u_s (\int_s^T u_t D_sG(t,T)dt)ds}{(\int_0^T u_t G(t,T)dt)^2
 }\right)\\
 D^W_u H&=&-\frac{\int_0^T u_s (\int_s^T u_t D_sG(t,T)dt)ds}{(\int_0^T u_t G(t,T)dt)^2
 }\delta^W(u)\\
 &&-S_T\left(\int_0^T u_t G(t,T)dt-\frac{\int_0^T u_s (\int_s^T u_t D_sG(t,T)dt)ds}{\int_0^T u_t G(t,T)dt
 }\right.\\
 &&\left.-\frac{\int_0^T u_r\int_r^T u_s (\int_s^T u_t D_rD_sG(t,T)dt)ds}{(\int_0^T u_t G(t,T)dt)^2}
 +2\frac{\left(\int_0^T u_s (\int_s^T u_t D_sG(t,T)dt)ds\right)^2}{(\int_0^T u_t G(t,T)dt)^3}\right)\\
 \partial_x H &=&-\frac{1}{x}S_T\left(1-\frac{\int_0^T u_s (\int_s^T u_t D_sG(t,T)dt)ds}{(\int_0^T u_t G(t,T)dt)^2
 }\right),
 \end{eqnarray*}
 $D^W_s G(t,T)$ and $D^W_rD^W_sG(t,T)$ are given by formulas $(\ref{derG1})$ and $(\ref{derG2})$.

\end{document}